\documentclass{amsart}
\usepackage{amsopn}
\usepackage{amssymb, amscd}

\newcommand{\nc}{\newcommand}

\nc{\fg}{\mathfrak{f} } \nc{\vg}{\mathfrak{v} } \nc{\wg}{\mathfrak{w} }
\nc{\zg}{\mathfrak{z} } \nc{\ngo}{\mathfrak{n} } \nc{\kg}{\mathfrak{k} }
\nc{\mg}{\mathfrak{m} } \nc{\bg}{\mathfrak{b} } \nc{\ggo}{\mathfrak{g} }
\nc{\ggob}{\overline{\mathfrak{g}} } \nc{\sog}{\mathfrak{so} }
\nc{\sug}{\mathfrak{su} } \nc{\spg}{\mathfrak{sp} } \nc{\slg}{\mathfrak{sl} }
\nc{\glg}{\mathfrak{gl} } \nc{\cg}{\mathfrak{c} } \nc{\rg}{\mathfrak{r} }
\nc{\hg}{\mathfrak{h} } \nc{\tg}{\mathfrak{t} } \nc{\ug}{\mathfrak{u} }
\nc{\dg}{\mathfrak{d} } \nc{\ag}{\mathfrak{a} } \nc{\pg}{\mathfrak{p} }
\nc{\sg}{\mathfrak{s} }

\nc{\pca}{\mathcal{P}} \nc{\nca}{\mathcal{N}} \nc{\lca}{\mathcal{L}}
\nc{\oca}{\mathcal{O}} \nc{\mca}{\mathcal{M}} \nc{\tca}{\mathcal{T}}
\nc{\aca}{\mathcal{A}} \nc{\cca}{\mathcal{C}} \nc{\gca}{\mathcal{G}}
\nc{\sca}{\mathcal{S}} \nc{\hca}{\mathcal{H}} \nc{\bca}{\mathcal{B}}
\nc{\dca}{\mathcal{D}} \nc{\val}{\operatorname{val}}

\nc{\diag}{\mathbf{d}}

\nc{\SO}{\mathrm{SO}} \nc{\Spe}{\mathrm{Sp}} \nc{\Sl}{\mathrm{SL}}\nc{\Si}{\mathrm{S}}
\nc{\SU}{\mathrm{SU}} \nc{\Or}{\mathrm{O}} \nc{\U}{\mathrm{U}} \nc{\Gl}{\mathrm{GL}}
\nc{\Se}{\mathrm{S}} \nc{\Cl}{\mathrm{Cl}} \nc{\Spein}{\mathrm{Spin}}
\nc{\Pin}{\mathrm{Pin}} \nc{\G}{\mathrm{GL}_n(\RR)} \nc{\g}{\mathfrak{gl}_n(\RR)}
\nc{\Kca}{\mathrm{K}}\nc{\Gr}{\mathrm{Gr}}

\nc{\RR}{{\Bbb R}} \nc{\HH}{{\Bbb H}} \nc{\CC}{{\Bbb C}} \nc{\ZZ}{{\Bbb Z}}
\nc{\FF}{{\Bbb F}} \nc{\NN}{{\Bbb N}} \nc{\QQ}{{\Bbb Q}} \nc{\PP}{{\Bbb P}}

\nc{\vp}{\varphi} \nc{\na}{\nabla}

\nc{\vs}{\vspace{.2cm}} \nc{\vsp}{\vspace{1cm}} \nc{\ip}{\langle\cdot,\cdot\rangle}
\nc{\ipp}{(\cdot,\cdot)} \nc{\la}{\langle} \nc{\ra}{\rangle} \nc{\unm}{\tfrac{1}{2}}
\nc{\unc}{\tfrac{1}{4}} \nc{\und}{\tfrac{1}{16}} \nc{\no}{\vs\noindent}
\nc{\lam}{\Lambda^2(\RR^n)^*\otimes\RR^n} \nc{\tangz}{{\rm T}^{\rm Zar}}
\nc{\nor}{{\sf n}}  \nc{\mum}{/\!\!/} \nc{\kir}{/\!\!/\!\!/}
\nc{\Ri}{\tfrac{4\Ric_{\mu}}{||\mu||^2}} \nc{\ds}{\displaystyle}
\nc{\ben}{\begin{enumerate}} \nc{\een}{\end{enumerate}} \nc{\f}{\frac}
\nc{\lb}{[\cdot,\cdot]} \nc{\isn}{\tfrac{1}{||v||^2}}

\nc{\Hess}{\operatorname{Hess}} \nc{\ad}{\operatorname{ad}}
\nc{\Ad}{\operatorname{Ad}} \nc{\rank}{\operatorname{rank}}
\nc{\Irr}{\operatorname{Irr}} \nc{\End}{\operatorname{End}}
\nc{\Aut}{\operatorname{Aut}} \nc{\Inn}{\operatorname{Inn}}
\nc{\Der}{\operatorname{Der}} \nc{\Ker}{\operatorname{Ker}}
\nc{\Iso}{\operatorname{I}} \nc{\Diff}{\operatorname{Diff}}
\nc{\Lie}{\operatorname{L}} \nc{\tr}{\operatorname{tr}} \nc{\dif}{\operatorname{d}}
\nc{\sen}{\operatorname{sen}} \nc{\modu}{\operatorname{mod}}\nc{\Id}{\operatorname{Id}}
\nc{\Ric}{\operatorname{R}} \nc{\Ricci}{\operatorname{Rc}}
\nc{\sym}{\operatorname{sym}} \nc{\symac}{\operatorname{sym^{ac}}}
\nc{\symc}{\operatorname{sym^{c}}} \nc{\scalar}{\operatorname{sc}}
\nc{\grad}{\operatorname{grad}} \nc{\ricci}{\operatorname{Ric}}
\nc{\Nor}{\operatorname{Norm}} \nc{\riccic}{\operatorname{ric^{c}}}
\nc{\riccig}{\operatorname{ric^{\gamma}}} \nc{\Rin}{\operatorname{M}}
\nc{\Le}{\operatorname{L}} \nc{\tang}{\operatorname{T}}
\nc{\level}{\operatorname{level}} \nc{\rad}{\operatorname{r}}
\nc{\abel}{\operatorname{ab}} \nc{\CH}{\operatorname{CH}}
\nc{\mcc}{\operatorname{mcc}} \nc{\Adj}{\operatorname{Adj}}
\nc{\Order}{\operatorname{O}}

\theoremstyle{plain}
\newtheorem{theorem}{Theorem}[section]
\newtheorem{proposition}[theorem]{Proposition}

\newtheorem{lemma}[theorem]{Lemma}

\theoremstyle{definition}
\newtheorem{definition}[theorem]{Definition}

\theoremstyle{remark}
\newtheorem{remark}[theorem]{Remark}

\title{The space of solsolitons in low dimensions}

\author{Cynthia Will}

\address{CIEM and FaMAF, Universidad Nacional de C\'ordoba, 5000 C\'ordoba, Argentina}
\email{cwill@famaf.unc.edu.ar}

\thanks{1991 {\it Mathematics Subject Classification.} Primary: 53C25, 53C30; Secondary: 22E25. \\
{\it Key words and phrases.} solitons, solvmanifold, left invariant metrics.\\
Partially supported by CONICET and research grants from FONCyT and SeCyT UNC
(Argentina).}

\begin{document}

\maketitle

\begin{abstract}
Up to now, the only known examples of homogeneous nontrivial Ricci soliton metrics are the so called {\it solsolitons}, i.e. certain left invariant metrics on simple connected solvable Lie groups. In this paper, we describe the moduli space of solsolitons of dimension $ \le 6$ up to isomorphism and scaling. We start with the already known classification of nilsolitons and, following the characterization given by Lauret in \cite{solsolitons}, we describe the subspace of solsolitons associated to a given nilsoliton, as the quotient of a Grassmanian by a finite group.
\end{abstract}

\section{Introduction}\label{intro}

A complete Riemannian metric $g$ on a differentiable manifold $M$ is said to be a {\it Ricci soliton} if its Ricci tensor satisfies
\begin{equation}
\Ricci(g)=c \,g +L_X\,g \qquad \mbox{for some}\; c\in \RR, \;\; X\in \chi(M)\, \mbox{complete},
\end{equation}
where $L_X$ denotes the Lie derivative in the direction of the vector field $X.$
Ricci soliton metrics came up in the study of the Ricci flow since they are the fixed points of the flow up to isometry and scaling. In the homogeneous case, all known nontrivial examples (i.e. not the product of an Einstein homogeneous manifold with a euclidean space) are isometric to a left-invariant metric $g$ on a simply connected solvable Lie group $G.$ Moreover, by identifying $g$ with an inner product on the Lie algebra $\ggo$ of $G,$ one has
\begin{equation}\label{eq1}
\ricci(g)=c I + D \qquad \mbox{for some} \, c\in \RR, \;\; D\in \Der(\ggo),
\end{equation}
where $\ricci(g)$  denotes the Ricci operator of $g.$ These metrics have been deeply studied, specially in the nilpotent case, where they are called {\it nilsolitons} (see for example \cite{finding},\cite{lw},\cite{N1},\cite{N2},\cite{P}).  In this case, Lauret proved many nice properties: they are the nilpotent parts of Einstein solvmanifolds, they are unique up to isometry and scaling and they are the critical points of certain natural variational problem (see the survey \cite{cruzchica}).

If $G$ is  solvable, metrics for which (\ref{eq1}) holds are called {\it solsolitons} and have been recently characterized by Lauret in \cite{solsolitons}.
 It is proved  there that any solsoliton is the extension of a nilsoliton by an abelian subalgebra of symmetric derivations of its metric Lie algebra. They are also unique in the sense that, among all left invariant metrics on a given simply connected solvable Lie group, there is at most one solsoliton up to isometry and scaling (see \cite{solsolitons}).

The aim of this paper is to approach the classification of solsolitons in low dimensions.  We first show in Section $3$ that, given an $n$-dimensional nilsoliton $N$, the space of solsolitons which are extensions of $N$ of dimension $r+n$ is parameterized, up to isometry and scaling, by
$$
\Gr_r(\ag)/W.
$$
Here $\ag$ is a maximal abelian subalgebra of symmetric derivations of the Lie algebra of $N,$  $\Gr_r(\ag)$ denotes the Grassmannian of
$r$-dimensional subspaces of $\ag$ and $W$ is a finite group. Thus for each $0 \le r \le \rank (N):= \dim(\ag),$ each nilsoliton $N$ provides a family of $(r+n)$-dimensional solsolitons depending on $r(\rank(N)-r)$ parameters. The number $\rank(N),$ called the {\it rank} of $N,$ is therefore the crucial datum to know on a nilsoliton if one is interested in classifying solsolitons.

{\small
\begin{table}
$$
\begin{array}{l|l|c}
\hline\hline
&&\\
{\rm nilsoliton\; metric} & {\rm eigenvalue \; type }& {\rm rank}  \\
&&\\
\hline &&\\
(0)  & (1;1) & 1\\
&&\\
(0,0)& (1;2)&2\\
&&\\
(0,0,0) & (1;3)& 3\\
&&\\
(0,0,\mathbf{12}) & (1<2;2,1) & 2 \\
 &&\\
\hline\hline
\end{array}
$$
\caption{Nilsolitons of dimension $\le 3$}\label{123dn}
\end{table} }

In Section $4$, we compute the rank of all nilsolitons of dimension $\le 6$ (see Tables 1-5). In most cases, we give in addition an explicit description of the coset $\Gr_r(\ag)/W$ as a real semialgebraic set (see \cite{Lafuente} for a similar approach in the case of nilsolitons attached to graphs). All this allows us to present quite a detailed picture of the moduli space $\mathcal{S}ol(m)$
of $m$-dimensional solsolitons for any $m \le 6$ (see Section $5$).
Finally, in Section $6,$ we make some remarks on negatively curved (sectional and Ricci) solsolitons.

{\bf Acknowledgments:} I would like to thank J. Lauret for his invaluable help and Comisi\'on Fulbright and CONICET for a three month fellowship when part of this research was done.

\vs

{\small
\begin{table}
$$
\begin{array}{c|l|l|c}
\hline\hline
&&&\\
N^o &{\rm nilsoliton\; metric \;\; } \eta_i &\quad {\rm eigenvalue \; type }& {\rm rank}  \\
&&&\\
\hline &&&\\
1 &(0,0,\mathbf{12} ,\mathbf{13} )  & (1<2<3<4;1,...,1)& 2\\
&&\\
2&(0,0,0,\mathbf{12}  ) & (2<3<4;2,1,1)& 3\\
&&&\\
3&(0,0,0,0) &
 (1;4)& 4 \\
 &&&\\
\hline\hline
\end{array}
$$
\caption{$4$-dimensional nilsolitons}\label{4dn}
\end{table} }

\section{Preliminaries}

Let $S$ be a solvmanifold, that is, a simply connected solvable Lie group endowed with a left invariant Riemannian metric. We will identify $S$ with its metric Lie algebra $(\sg \,, \ip),$ where $\sg$ is the Lie algebra of $S$ and $ \ip $ denotes the inner product on $\sg$ determined by the metric. If $\ngo$ is the nilradical of $\sg$ (i.e. the maximal nilpotent ideal), consider the orthogonal decomposition
$$ \sg=\ag \oplus \ngo.$$
$S$ is called {\it standard} if $[\ag,\ag]=0.$
\begin{definition}\label{defsols}
A left-invariant metric $g$ on $S$ is said to be a
{\it solsoliton} if
\begin{equation}\label{rsD}
\ricci(g)=cI+D, \qquad\mbox{for some}\quad c\in\RR,\quad D\in\Der(\sg),
\end{equation}
where $\ricci(g)$ denotes the {\it Ricci operator} of $g.$
\end{definition}
Hence solsolitons are natural generalizations of Einstein metrics (i.e. when $D=0$). Einstein solvmanifolds have been studied by Heber who gave in \cite{H} many structural results. In particular, if $S$ is an Einstein solvmanifold, he showed that there exists an element $H\in \ag$ such that the eigenvalues of $\ad H_{|\ngo}$ are positive integers without common divisors, say $k_1<\dots <k_r$ with multiplicities $d_1,\dots,d_r.$ The {\it eigenvalue type} of $S$ is defined as the tuple
 $$ (k;d)=(k_1< \dots <k_r ; d_1,\dots,d_r). $$
We have that $\ad H_{|\ngo}$ is a multiple of $D$ in equation (\ref{rsD}).

When the group $S$ is nilpotent, metrics for which (\ref{rsD}) holds are called
{\it nilsolitons}. They have been deeply studied by Lauret (see \cite{cruzchica}) who has proved, among others properties, that:
 \ben \item[$\bullet$] a given nilpotent Lie group $N$ can admit at most one nilsoliton up to isometry and scaling among all its left-invariant metrics,
 \item[$\bullet$]  nilsolitons are precisely the metric nilradicals of Einstein solvmanifolds,
\item[$\bullet$] any Ricci soliton left invariant metric is a nilsoliton.
\een
By this uniqueness result we can associate to a nilsoliton the eigenvalue type of its uniquely defined rank-one Einstein solvable extension.

We give in what follows a summary of Lauret's results we are using in our classification on solsolitons (see \cite{solsolitons}). Note that we have slightly changed the notation in \cite{solsolitons} to be consistent with the one that is useful in this context.
For any metric nilpotent Lie algebra $(\ngo, \ip),$ let $\sym(n)$ denote the space of symmetric transformations of $\ngo$ with respect to $\ip$.

It is proved in \cite{solsolitons} that any solsoliton can be constructed, up to isometry and scaling, from a nilsoliton $(\ngo,\ip)$ and an abelian subalgebra of symmetric derivations of $\ngo$ in the following way (see \cite[Corollary 4.10]{solsolitons}).

\begin{proposition}\label{const}\cite[Proposition 4.3]{solsolitons}
Let $(\ngo,\ip_1)$ be a nilsoliton with Ricci operator $\ricci_1=cI+D_1$,
$c<0$, $D_1\in\Der(\ngo)$, and consider $\ag$ any abelian Lie algebra of symmetric
derivations of $(\ngo,\ip_1)$.  Then the solvmanifold $S$ with Lie algebra
$\sg=\ag\oplus\ngo$ (semidirect product) and inner product given by
$$
\ip|_{\ngo\times\ngo}=\ip_1, \qquad \la\ag,\ngo\ra=0, \qquad \la
D,D\ra=-\tfrac{1}{c}\tr{D^2} \qquad\forall\, D\in\ag,
$$
is a solsoliton with $\ricci=cI+D_0$, where $D_0\in\Der(\sg)$ is defined by
$D_0|_{\ag}=0$, $D_0|_{\ngo}=D_1-\ad{H}|_{\ngo}$ and $H$ is the mean curvature vector of
$S$. Furthermore, $S$ is Einstein if and only if $D_1\in\ag$.
\end{proposition}

Moreover, as in the case of nilsolitons, a given solvable Lie group admits at most one Ricci soliton left invariant metric up to isometry and scaling (see \cite[Theorem 5.1]{solsolitons}. Therefore, to classify solsolitons we need to have first a classification of nilsolitons and for each one of them we have to study the subalgebras of symmetric derivations up to the action of the orthogonal automorphism group, since the equivalence among them is given by the following proposition.

\begin{proposition}\label{toclas}\cite[Proposition 5.3]{solsolitons}
Let $(\ngo,\ip)$ be a nilsoliton and let $S$ and $S'$ be two solsolitons
constructed as in {\rm Proposition \ref{const}} for abelian subalgebras
$$
\ag,\ag'\subset\Der(\ngo)\cap\sym(n),
$$
respectively.  Then $S$ is isometric to $S'$ if and only if there exists
$h\in\Aut(\ngo)\cap\Or(\ip)$ such that $\ag'=h\ag h^{-1}$.
\end{proposition}

{\small
\begin{table}
$$
\begin{array}{c|l|l|c}
\hline\hline &&& \\
N^o &\qquad {\rm nilsoliton\; metric \;\; } \lambda_i &\qquad {\rm eigenvalue \; type }&{\rm rank} \\ &&& \\
\hline 
1 &\Big(0,0,3\, \mathbf{12} ,4\,\mathbf{13} ,3\,\mathbf{14} \Big)  & (2<9<11<13<15;1,...,1)& 2\\

2&\Big(0,0,3^{1/2} \mathbf{12}, 3^{1/2} \mathbf{13}, 2^{1/2} \mathbf{14} + 2^{1/2} \mathbf{23}\Big) &  (1<2<3<4<5;1,...,1) &1\\

3&\Big(0,0,0, 2^{1/2}\mathbf{12} ,\mathbf{23} +  2^{1/2}\mathbf{14} \Big) & (3<4<6<7<10;1,...1) &2\\

4&\Big(0,0,0,0,\mathbf{12} +\mathbf{34} \Big)
 & (1<2;4,1) &3 \\

5&\Big(0,0,2\,{\mathbf  12},3^{1/2}\mathbf{13} ,3^{1/2} \mathbf{23}\Big) &
(1<2<3;2,1,2)&2 \\

6&\Big(0,0,0,\mathbf{12} ,\mathbf{13} \Big) & (2<3<5;1,2,2)&3 \\

7&\Big(0,0,0,0,\mathbf{12} \Big) & (2<3<4;2,2,1)&4 \\

8&\Big(0,0,0, \mathbf{12} , \mathbf{14} \Big) & (1<2<3<4;1,1,2,1)&3 \\

9&\Big(0,0,0,0,0\Big) & (1;5)&5 \\

\hline\hline
\end{array}
$$
\caption{$5$-dimensional nilsolitons}\label{5dn}
\end{table} }

Let $\mu: \RR^n \times \RR^n \to \RR^n$ be a skew-symmetric bilinear form.  If $\mu$ satisfies the Jacobi identity then it defines a Lie algebra $(\RR^n, \mu)$. We will denote by
  $$
\nca=\{\,\mu:\RR^n \times \RR^n \to \RR^n:\mu\;\mbox{skew-symmetric, bilinear, nilpotent and satisfies Jacobi}\}.
$$
$\nca$ is often called the variety of nilpotent Lie algebras (of dimension $n$).
If we also fix the canonical basis  on $\RR^n,$ $\{X_1, \dots, X_n\},$ and the inner product $\ip$ that makes this basis orthonormal then each $\mu \in \nca$ define a metric nilpotent Lie algebra $(\RR^n, \mu, \ip).$ We will say then that $\mu \in \nca$ is a nilsoliton if $(\RR^n, \mu, \ip)$ is so.

 Nilsolitons metrics are classified up to dimension $6$ (see \cite{finding} and \cite{W}). They are listed in Tables \ref{123dn} through \ref{pc23}.
 In these tables, an element $\mu \in \nca$ is represented as a vector of $n$ coordinates in such a way that
$$\begin{array}{lcl}
k^{\mbox{th}}\mbox{-coordinate}\,=\, c\,\mathbf{ij}+d\,\mathbf{rl}& \Leftrightarrow & \mu(X_i,X_j)=c X_k,\,\, \mu(X_r,X_l)=d X_k.
\end{array}$$
Hence, for example, $(0,0,0,2^{\unm}\mathbf{12}, \mathbf{23}+2^{\unm}\mathbf{14})$ represents the  $5$-dimensional Lie algebra with Lie bracket given by
$$
\mu(X_1,X_2)=2^{\unm} X_4, \qquad  \mu(X_2,X_3)=X_5, \qquad  \mu(X_1,X_4)=2^{\unm}X_5.
$$

We will denote by
$\diag(a_1,\dots,a_n)$ the diagonal matrix with entries $a_1, \dots, a_n,$ and more in general
\begin{equation}\label{diagonal}
\diag(M_1,\dots, M_r):= \left[\begin{smallmatrix}
  M_1&&&\\&\ddots&\\&&M_r \end{smallmatrix}\right],
  \end{equation}
  for matrices $M_i \in \glg(n_i)$ such that $n_1+\dots+n_r=n.$

{\small
\begin{table}
$$
\begin{array}{c|l|l|c}
\hline\hline &&& \\
N^o &\qquad {\rm Nilsoliton \; metric: \;\; } \mu_i &\qquad {\rm eigenvalue \; type }&{\rm rank }\\ && \\
\hline 
1 &\Big(0,0,(13)^\frac{1}{2}\mathbf{12},4\,\mathbf{13},(12)^\frac{1}{2}\mathbf{14}+2\,\mathbf{23},  & (1<2<3<4<5<7;1,...,1)&1\\
 &\qquad (12)^\frac{1}{2}\mathbf{34}+(13)^\frac{1}{2}\mathbf{52}\Big)& &  \\

2&\Big(0,0,\mathbf{12},(\frac{4}{3})^\frac{1}{2}\mathbf{13},\mathbf{14},\mathbf{34}+\mathbf{52}\Big) & (1<3<4<5<6<9;1,...,1)&2 \\

3&\Big(0,0,2\,\mathbf{12},6^\frac{1}{2}\mathbf{13},6^\frac{1}{2}\mathbf{14},2\, \mathbf{15}\Big) & (1<9<10<11<12<13;1,...1)&2 \\

4&\Big(0,0,(22)^\frac{1}{2}\mathbf{12},6\;\mathbf{13},(22)^\frac{1}{2}\mathbf{14}+(30)^\frac{1}{2}\mathbf{23},
 & (1<2<3<4<5<6;1,...1) &1 \\
 & \qquad 5\,\mathbf{24}+(30)^\frac{1}{2}\mathbf{15}\Big) && \\

5&\Big(0,0,7^\frac{1}{2}\mathbf{12},(\frac{15}{2})^\frac{1}{2}\mathbf{13},3\;\mathbf{14}, &
(1<3<4<5<6<7;1,...1)&1 \\
 &\qquad (\frac{15}{2})^\frac{1}{2}\mathbf{23}+2\,\mathbf{15}\Big)&& \\

6&\Big(0,0,\mathbf{12},\mathbf{13},\mathbf{23},\mathbf{14}\Big) & (1<2<3<4<5;1,1,1,1,2) &2\\

7&\Big(0,0,2\;\mathbf{12}, 5^\frac{1}{2}\mathbf{13},5^\frac{1}{2}\mathbf{23},2\;\mathbf{14}-2\;\mathbf{25}\Big) & (1<2<3<4;2,1,2,1) &2\\

8&\Big(0,0,2\;\mathbf{12}, 5^\frac{1}{2}\mathbf{13},5^\frac{1}{2}\mathbf{23},2\;\mathbf{14}+2\;\mathbf{25}\Big) & (1<2<3<4;2,1,2,1)&1 \\

9&\Big(0,0,0,(\frac{5}{4})^\frac{1}{2}\mathbf{12},\mathbf{14}-\mathbf{23},(\frac{5}{4})^\frac{1}{2}\mathbf{15}+\mathbf{34}\Big) &
 (6<11<12<17<23<29;1,...,1)&2 \\

10&\Big(0,0,0,\mathbf{12},(\frac{5}{3})^\frac{1}{2}\mathbf{14},\mathbf{15}+\mathbf{23}\Big) & (4<9<12<13<17<21;1,...,1)&2\\

11&\Big(0,0,-(\frac{35}{17\cdot 8})^\frac{1}{2}\mathbf{12},(\frac{21}{17\cdot 2})^\frac{1}{2}\mathbf{12},
(\frac{25}{17\cdot 4})^\frac{1}{2}\mathbf{14} & (1<2<3<4<5;1,1,2,1,1)&1 \\

 &\qquad - (\frac{15}{17})^\frac{1}{2} \mathbf{13},(\frac{3}{4})^\frac{1}{2}\mathbf{15}+(\frac{7}{8})^\frac{1}{2}\mathbf{24}\Big) &   \\

12&\Big(0,0,0,3^\frac{1}{2}\mathbf{12},3^\frac{1}{2}\mathbf{14},2^\frac{1}{2}\mathbf{15}+2^\frac{1}{2}\mathbf{24}\Big) &
(3<6<9<11<12<15;1,...,1) &2\\

13&\Big(0,0,0,3^\frac{1}{2}\mathbf{12},2\;\mathbf{14},3^\frac{1}{2}\mathbf{15}\Big) & (2<9<11<12<13<15;1,...,1)&3\\

\hline\hline
\end{array}
$$
\caption{$4$ and $5$-step $6$-dimensional nilsolitons}\label{pc45}
\end{table} }

{\small
 \begin{table}
 $$
\begin{array}{c|l|l|c}
\hline\hline  &&\\
N^o &\qquad {\rm Nilsoliton \; metric:\;\;\;} \mu_i &\qquad {\rm eigenvalue \;\; type } &{\rm rank} \\&&&\\
 \hline 

 14 &\Big(0,0,0,3^\frac{1}{2}\mathbf{12},2^\frac{1}{2}\mathbf{13},2^\frac{1}{2}\mathbf{14}+3^\frac{1}{2}\mathbf{35}\Big) &
 (2<3<4<5<6<8;1,...,1)&2 \\

15& \Big(0,0,0,\mathbf{12},\mathbf{23},\mathbf{14}+\mathbf{35}\Big) & (1<2<3;3,2,1)& 3 \\

16&\Big(0,0,0,\mathbf{12},\mathbf{23},\mathbf{14}-\mathbf{35}\Big) &  (1<2<3;3,2,1)&2 \\

17&\Big(0,0,0,2\, \mathbf{12},3^\frac{1}{2}\mathbf{14},3^\frac{1}{2}\mathbf{24}\Big) & (5<10<12<15;2,1,1,2)&3  \\

18&\Big(0,0,0,2^\frac{1}{2}\mathbf{12},2^{-\frac{1}{2}}\mathbf{13}+ (\frac{3}{2})^\frac{1}{2}\mathbf{42}, & (1<2<3;2,2,2) &1\\
& \qquad (\frac{3}{2})^\frac{1}{2}\mathbf{14}+ 2^{-\frac{1}{2}}\mathbf{23}\Big) && \\

19 &\Big(0,0,0,2\, \mathbf{12},3^\frac{1}{2}\mathbf{14},13+3^\frac{1}{2}\mathbf{42}\Big) & (5<6<11<12<16<17;1,...,1)
&2\\

20&\Big(0,0,-\mathbf{12},3^\frac{1}{2}\mathbf{12},2\,\mathbf{14}, \mathbf{24}-3^\frac{1}{2}\mathbf{23}\Big) & (1<2<3;2,2,2) &2\\

21&\Big(0,0,0,2^\frac{1}{2}\mathbf{12},\mathbf{13},2^\frac{1}{2}\mathbf{14}+\mathbf{23}\Big) & (3<5<6<8<9<11;1,...,1)&2\\

22&\Big(0,0,0,(\frac{3}{4})^\frac{1}{2}\mathbf{12},(\frac{3}{4})^\frac{1}{2}\mathbf{13},\mathbf{24}\Big) & (5<6<9<11<15<16;1,...,1)&3 \\

23&\Big(0,0,0,2^\frac{1}{2}\mathbf{12},\mathbf{13},2^\frac{1}{2}\mathbf{14}\Big) & (2<5<6<7<8<9;1,...,1)&3 \\

24&\Big(0,0,0,\mathbf{12},\mathbf{13},\mathbf{23}\Big) & (1<2;3,3)&3 \\

25&\Big(0,0,0,0,2\mathbf{12},3^\frac{1}{2}\mathbf{15}+ 3^\frac{1}{2}\mathbf{34}\Big) & (5<8<9<13<18;1,1,2,1,1)&3 \\

26&\Big(0,0,0,0,\mathbf{12},\mathbf{15}\Big) & (1<2<3<4;1,1,3,1)&4 \\

27&\Big(0,0,0,0,2^\frac{1}{2}\mathbf{12},\mathbf{14}+2^\frac{1}{2}\mathbf{25}\Big) & (3<4<6<7<10;1,1,1,2,1)&3 \\

28&\Big(0,0,0,0,\mathbf{13}+\mathbf{42},\mathbf{14}+\mathbf{23}\Big) & (1<2;4,2) &2\\

29&\Big(0,0,0,0,\mathbf{12},\mathbf{14}+\mathbf{23}\Big) & (3<4<6<7;2,2,1,1) &3\\

30&\Big(0,0,0,0,\mathbf{12},\mathbf{34}\Big) & (1<2;4,2)&4 \\

31&\Big(0,0,0,0,\mathbf{12},\mathbf{13}\Big)& (2<3<4<5;1,2,1,2) &4 \\

32&\Big(0,0,0,0,0,\mathbf{12}+\mathbf{34}\Big) & (3<4<6;4,1,1)&4 \\

33&\Big(0,0,0,0,0,\mathbf{12}\Big) & (2<3<4;2,3,1)&5 \\

34&\Big(0,0,0,0,0,0) & (1;6)&6 \\

 \hline\hline
\end{array}
$$
\caption{$k$-step $6$-dimensional nilsolitons, $k\le 3$.}\label{pc23}
\end{table} }

\section{General facts on the classification of solsolitons}

Let $(\RR^n,\mu)$ be a nilpotent Lie algebra. We fix the canonical inner product $\ip$ on $\RR^n$ and denote by $\sog(n),\; \sym(n)$ and $\Or(n)$ the subspaces of skew-symmetric, symmetric linear maps on $\RR^n$ and the orthogonal group relative to $\ip,$ respectively. Let $G_\mu$ denote the closed subgroup of automorphisms given by
$$
G_\mu:=\{g \in \Aut(\mu): g^t \in \Aut(\mu)\}.
$$
Since $\Aut(\mu)$ is algebraic, we have that $G_\mu$ is a real reductive group (see \cite[Section 2.1]{Wa}) with Lie algebra
$$
\ggo_\mu=\{A \in \Der(\mu): A^t \in \Der(\mu)\},
$$
and Cartan decompositions
$$ G_\mu = \Kca_\mu \exp(\pg_\mu), \qquad \ggo_\mu=\kg_\mu \oplus \pg_\mu,$$
\noindent where
$$
\begin{array}{c}
\Kca_\mu:= G_\mu \cap \Or(n)=\Aut(\mu)\cap \Or(n), \qquad
\kg_\mu:= \ggo_\mu \cap \sog(n)=\Der(\mu)\cap \sog(n),\\ \\
\pg_\mu:= \ggo_\mu \cap \sym(n)=\Der(\mu)\cap \sym(n).
\end{array}
$$

 It is proved in \cite[Lemma 4.7]{solsolitons}) that $\ggo_\mu$ is precisely the space of derivations of $(\RR^n,\mu)$ which are normal (i.e. $[D,D^t]=0$).

Let us assume from now on that $(\RR^n,\mu,\ip)$ is a nilsoliton.
According to Proposition \ref{toclas}, the solsolitons associated with the nilsoliton $\mu$ are parameterized by abelian subspaces of $\pg_\mu$ up to conjugancy by $\Kca_\mu.$ Let $\ag_\mu$ be a maximal abelian subalgebra of $\pg_\mu.$ Since $G_\mu$ is real reductive, it is well known that any abelian subspace of $\pg_\mu$
is $\Kca_\mu$-conjugate to a subspace of $\ag_\mu$ (see \cite[2.1.9]{Wa}). Moreover, any two subspaces of $\ag_\mu$ are  $\Kca_\mu$-conjugate if and only if they are conjugate by an element of the Weyl group $W_\mu$ of $G_\mu$ (see \cite[Proposition 2.2, Ch.VII]{Hel}), which is given by
$$
W_\mu:=N_{\Kca_\mu}(\ag_\mu)/Z_{\Kca_\mu}(\ag_\mu),
$$
\noindent where
$$\begin{array}{c}
N_{\Kca_\mu}(\ag_\mu):=\{g \in \Kca_\mu: \Ad(g)\ag_\mu \subset \ag_\mu\}, \\ \\
Z_{\Kca_\mu}(\ag_\mu):=\{g \in \Kca_\mu: \Ad(g)A=A,\;\;\; \forall A \in \ag_\mu\}.
\end{array}
$$

We therefore obtain that the set of $(r+n)$-dimensional solsolitons associated with the nilsoliton $\mu$ is parameterized, up to isometry and scaling, by the coset
\begin{equation}\label{para}
\Gr_r(\ag_\mu)/W_\mu, 
\end{equation}
\noindent where $\Gr_r(\ag_\mu)$ is the Grassmanian of  $r$-dimensional subspaces of $\ag_\mu$. We note that $0 \le r \le \rank(\mu),$ where $\rank(\mu):= \dim \ag_\mu$ will be called the {\it rank} of the nilsoliton $\mu$. As $W_\mu$ is a finite group, the quotient $\Gr_r(\ag_\mu)/W_\mu$ depends on $r(\rank(\mu)-r$)
parameters.

The aim of this paper is to compute the rank for  any nilsoliton of dimension $\le 6,$ and also to give in most cases an explicit description of the quotient in (\ref{para}).

We finish this section by dealing with the cases that admit a simultaneous study. In the next section we will work out a case by case analysis of the remaining nilsolitons.

\vspace{.3cm}

\subsection{Abelian nilsolitons}

A first case that can be studied in general is the abelian case, that is,  $\ngo=(\RR^n, \mu)$ with $\mu=0.$ We have that $\Der(0)= \glg(n),$ and therefore a maximal abelian subspace of symmetric derivations  can be chosen to be the set of diagonal matrices:
$$\pg_0=\sym(n) \qquad \ag_0=\{\diag(a_1,\dots,a_n): a_i \in \RR\}, \qquad \rank(0)=n.$$
 Since $\Kca_0 = \Or(n),$ the action by conjugation on diagonal matrices contains all the permutations of the entries. Hence, for $0\le r \le n,$ the space of $(r+n)$-dimensional solsolitons is given by
$$\Gr_r(\RR^n)/\Si_n,$$
\noindent where $\Si_n$ is the permutation group of $n$ elements. Note that the $(n+1)$-dimensional hyperbolic space $\RR H^{n+1}$ is an element of
$\Gr_r(\RR^n)/\Si_n = \PP\RR^n/\Si_n$ and is characterized as the only point that corresponds to an Einstein metric.

\subsection{Maximal and minimal dimension solsolitons}\label{maxrank}

If $(\RR^n, \mu, \ip)$ is a nilsoliton such that $\rank (\mu) = k,$ then
$$\Gr_0(\ag_{\mu})=\{*\}, \qquad \qquad \Gr_k(\ag_{\mu})=\{*\}.$$
In other words, there is only one solsoliton $\sg = \ag \oplus \ngo_\mu$ of dimension $(k+n)$ associated to the nilsoliton $\mu,$ which is always Einstein (see Proposition \ref{const}), and only one of dimension $n,$ the nilsoliton itself, which is never Einstein unless $\mu=0$.

 \subsection{Rank-one nilsolitons}\label{rank1}

 Another case we can study separately is when $\ag_\mu$ is one dimensional. Then
  $$\Gr_1(\ag_{\mu})=\{*\},$$
  that is, there is only one solsoliton associated to $\mu$ which is always Einstein. For dimension $\le 6$ they are all listed in Table \ref{sol6D1}.

{\small
 \begin{table}
 $$
\begin{array}{c|c|c}
\hline\hline  &\\
{\rm metric} & \ag_\mu & {\rm dim.} \\& &\\
 \hline &&\\
 \lambda_2 & \RR \diag(1,2,3,4,5)\;\; &5 \\
\mu_1 & \RR \diag(1,2,3,4,5,7)&6 \\
\mu_4 & \RR \diag(1,2,3,4,5,6)&6 \\
\mu_5 & \RR \diag(1,3,4,5,6,7)&6 \\
\mu_8 & \RR \diag(1,1,2,3,3,4)&6 \\
\mu_{11} & \RR \diag(1,2,3,3,4,5))&6 \\
\mu_{18}& \RR \diag(1,1,2,2,3,3))&6 \\
&\\
 \hline\hline
\end{array}
$$
\caption{Rank-one nilsolitons of dimension $\le 6$.}\label{sol6D1}
\end{table} }

\subsection{Multiplicity-free eigenvalue type nilsoliton}

It is well known that skew-symmetric and symmetric derivations, as well as orthogonal automorphisms of $(\RR^n,\mu,\ip),$ all commute with $\ricci_\mu$ (see \cite[Lemma 2.2]{H}). In this way, if $\mu$ is an nilsoliton, say with $\ricci_\mu=c \Id +D_\mu,$ all the above mentioned operators also commute with $D_\mu.$ Then if we fix $\beta=\{X_1, \dots X_n\}$ a basis of eigenvectors of $D_\mu$ we get the following observation.

\begin{lemma}\label{util}
Let $(\RR^n,\mu,\ip)$ be a nilsoliton with eigenvalue type $(k_1 < \dots < k_r; n_1, \dots,n_r).$ Then the elements in $\pg_\mu$ and in $\Kca_\mu$ are diagonal block matrices of dimension $n_1, \dots,n_r,$ as in {\rm (\ref{diagonal})}.
\end{lemma}

\begin{definition}\label{fet} An eigenvalue type of the form $(k_1< \dots < k_r; 1, \dots,1)$ will be called {\it multiplicity-free}.
\end{definition}

 If $(\RR^n,\mu,\ip)$ is a nilsoliton with multiplicity-free eigenvalue type, then the elements in $\pg_\mu$ are diagonal matrices and $\Kca_\mu \subseteq \diag(\pm 1,\dots, \pm 1).$ This implies that $\Kca_\mu$ acts trivially by conjugation on $\ag_\mu$ and consequently,
 the solsolitons of dimension $r+n$ are just parameterized by
 $$\Gr_r(\ag_\mu), \qquad 0\le r \le \rank(\mu).$$
 The multiplicity-free eigenvalue type nilsolitons of rank
 $\ge 2$ are listed in Table \ref{sol111}, where the relevant information is given.

{\small
 \begin{table}
 $$
\begin{array}{c|l|l|c|c}
\hline\hline  &&&&\\
{\rm metric} &\qquad \ag_\mu & {\rm Einstein \;condition } & {\rm rank } & {\rm dim.} \\&&&&\\
 \hline 

 \eta_1 &(a,b,a+b,2a+b)&2a=b&2& 4 \\

\lambda_1& (a,b,a+b,2a+b,3a+b)&  9a={2}b&2& 5\\


\lambda_3& (a,b,2a,a+b,2a+b)& 4a={3}b & 2&5 \\


\mu_2 & (a,b,a+b,2a+b,3a+b,3a+2b)& 3a=b & 2& 6\\

\mu_3& (a,b,a+b,2a+b,3a+b,4a+b)& 9a=b & 2& 6 \\



\mu_9 &(a,b,2a,a+b,2a+b,3a+b) &11 a=6 b& 2&6\\

\mu_{10}& (a,b,3a,a+b,2a+b,3a+b)& 9a={4}b&2&6\\

\mu_{12}& (a,2a,b,3a,a+b,2a+b)&11 a=3 b&2 & 6 \\

\mu_{13}&  (a,b,c,a+b,2a+b,3a+b)&a/2=b/9=c/12& 3& 6 \\

\mu_{14} & (a,b,\frac{a+b}{2},a+b,\frac{3a+b}{2},2a+b) & 2a=b & 2&6 \\

\mu_{19} &(a,b,2b,a+b,2a+b,a+2b) & 6 a=5 b& 2 &6 \\

\mu_{21}& (a,b,2a,a+b,3a,2a+b )& 5 a= 3 b& 2&6\\

\mu_{22}&(a,b,c,a+b,a+c,a+2b ) & a/6=b/5=c/9& 3&6 \\

\mu_{23}&( a,b,c,a+b,a+c,2a+b) & a/2=b/5=c/6& 3&6 \\

 \hline\hline
\end{array}
$$
\caption{Multiplicity-free eigenvalue type nilsolitons of rank $\ge 2$.}\label{sol111}
\end{table} }

\subsection{Rank $3$}\label{rank3}

It is easy to see that $\Gr_r(\RR^k) \simeq \Gr_{k-r}(\RR^k),$ for example, by assigning to an $r$-dimensional subspace of $\RR^k$ its orthogonal complement. Therefore, for a nilsoliton of rank $k$  the space of associated  $(r+n)$-dimensional solsolitons is isomorphic to the space of $(k-r+n)$-dimensional ones, that is,
$$
\Gr_r(\ag_\mu)/W_\mu \simeq \Gr_{k-r}(\ag_\mu)/W_\mu.
$$
In particular, if the rank of $(\RR^n,\mu)$ is $3,$ we only have to describe the space of $(n+1)$-dimensional solsolitons. In fact, there is only one $(n+3)$-dimensional solsoliton, the Einstein one, and the space of $(n+2)$-dimensional solsolitons is isomorphic, by the above observation, to the space of $(n+1)$-dimensional solsolitons.  Moreover, we have that
$$\Gr_1(\ag_\mu)\simeq \PP \ag_\mu \simeq \PP \RR^k,$$
\noindent where $\PP V$ denotes the projective space of $V$.

\section{Case by case analysis}

\vs
In this section, we study the nilsolitons of dimension $\le 6$ that are not included in Tables \ref{sol6D1} or \ref{sol111}. For each one of these nilsolitons $\mu,$ we give a set that parameterizes the space of $(n+1)$-dimensional solsolitons with nilradical $\mu$ (up to isometry and scaling) by explicitly describing $\PP\ag_\mu/W_\mu$ (see (\ref{para})).

Note that the nilsolitons we are studying have rank at most $5$ and therefore it follows from what was explained in Sections \ref{maxrank}, \ref{rank1} and \ref{rank3} that this is the only case we need to consider for those of rank one, two or three. If the rank is greater than three, to have a complete parametrization of the set of solsolitons associated to $\mu,$ it remains to consider $\Gr_2(\ag_\mu)/W_\mu,$ but this is a problem that exceeds the scope of this paper. We hope to
deal with it in a future work.

We will denote by $\ngo_\mu$ the Lie algebra $(\RR^n, \mu,\ip),$ where $\ip$ is the canonical inner product we have fixed,
and by $\sigma_{ij}$ the element of the symmetric group $S_n$ that permutes $X_i$ with $X_j$ and fixes all other elements of the basis.

\subsection{Nilsolitons of dimension $3$}\
This case has been considered in \cite{solsolitons} but we will include it again, for completeness.

$\bullet$ $\hg_3$, the $3$-dimensional Heisenberg Lie algebra, given by $(0,0,\mathbf{12})$ (see Table \ref{123dn}). Therefore
  $$\pg_{\hg_3}=\left\{ \left[\begin{smallmatrix} a&c&\\c&b&\\&&a+b \end{smallmatrix}\right]:\, a,b,c \in \RR \right\}, \qquad \Kca_{\hg_3}=\left\{ \diag(H,\det H)\,:\, H \in \Or(2) \right\},$$
   $$\ag_{\hg_3}=\left\{\diag(a,b,a+b)\,:\, a,b \in \RR \right\}.$$
   $W_{\hg_3}$ is the group generated by $\{\Id, \sigma_{12}\}$.

  \noindent The corresponding set $\PP\ag_{\hg_3}/ W_{\hg_3}$ of $4$-dimensional solsolitons is parameterized by:
 \begin{equation}\label{f1}
\mathfrak{F}^1_1= \left\{(a,b) \in S^1:\; \begin{array}{ll}
  |a|\leq b
\end{array}\right\},
\end{equation}
where, in general, $S^{n}$ is the $n$-dimensional sphere in $\RR^{n+1}.$

\noindent Einstein condition: $a=b.$

\noindent This correspondence assigns to each $(a,b) \in \mathfrak{F}^1_1$ the solsoliton  $\sg = \RR \left[
 \begin{smallmatrix} a&& \\&b&\\
&&a+b\end{smallmatrix}\right] \oplus  \hg_3,$ where Lie bracket and inner product are respectinely given by
$$[X_1,X_2]=X_3,\;\; [X_0,X_1]=a\;X_1,\;[X_0,X_2]=b\,X_2,\;\;[X_0,X_3]=(a+b)\;X_3,$$
$$\langle X_i,X_j \rangle= \delta_{ij},\;\;\; \mbox{for}\,\; 0\le i,j \le 3,\;\; \langle X_0,X_0 \rangle = \frac{4}{3}(a^2+b^2+ab). \hspace{1.5cm}$$

\subsection{Nilsolitons of dimension $4$}(see Table \ref{4dn})

$\bullet$ $\mathbf{\eta}_2$  ($\simeq \hg_3  \oplus \RR X_3$) Derivations of $\eta_2$ have been studied in \cite{dim3y4}.
$$\pg_{\eta_2}=\left\{ \diag(A,c,\tr A))\,:\, A \in \sym(2), c \in \RR \right\}, \quad \Kca_{\eta_2}=\left\{ \diag(H,\pm 1,\det H)\,:\, H \in \Or(2) \right\},$$
  $$\ag_{\eta_2}=\left\{\diag(a,b,c,a+b)\,:\, a,b,c \in \RR \right\}. $$
  $W_{\eta_2}$ is the group generated by $\{\Id, \sigma_{12}\}$.

  \noindent The corresponding $5$-dimensional solsolitons are parameterized by:
\begin{equation}\label{f2}
\mathfrak{F}^2_1= \left\{(a,b,c) \in S^2:\; \begin{array}{ll}
  c>0 &  a\leq b, \, \mbox{or}\\ c=0 & |a|\leq b
\end{array}\right\}.
\end{equation}
Einstein condition: $a=b=2/3 c.$

\noindent The above correspondence is given by
$$(a,b,c) \mapsto \RR \left[
 \begin{smallmatrix} a&&& \\&b&&\\
&&c&\\ &&&a+b\end{smallmatrix}\right] \oplus \ngo,$$
where $\ngo=(\RR^4, \eta_2).$

The Lie bracket and the inner product of $\sg$ are explicitly given by:
$$[X_1,X_2]=X_4,\,[X_0,X_1]=a\,X_1,\,[X_0,X_2]=b\,X_2,\,[X_0,X_3]=c\,X_3,\,[X_0,X_4]=(a+b)\,X_4,$$
$$\langle X_i,X_j \rangle= \delta_{ij},\;\;\; \mbox{for}\,\; 0\le i,j \le 4,\;\; \langle X_0,X_0 \rangle = \frac{4}{3}(a^2+b^2+ab+\frac{c^2}{2}). \hspace{3.5cm}$$

\subsection{Nilsolitons of dimension $5$}(see Table \ref{5dn})\

\vs

$\bullet$ $\mathbf{\lambda}_4$ ($5$-dimensional Heisenberg Lie algebra).
$$\begin{array}{c}
\pg_{\lambda_4}=\left\{\left[\begin{smallmatrix}\begin{smallmatrix}\begin{smallmatrix}a&\alpha\\ \alpha & b \end{smallmatrix} & B\\ B &
\begin{smallmatrix} c &\beta\\ \beta & a+b-c\end{smallmatrix} \end{smallmatrix}& \\ &a+b \end{smallmatrix} \right]
\,:\, a,b,c,\alpha,\beta \in \RR, \; B= \left[\begin{smallmatrix}u&v\\v & -u \end{smallmatrix} \right] \right\},\\
\ag_{\lambda_4} = \left\{ \diag(a,b,c,a+b-c,a+b)\,:\, a,b,c \in \RR \right\}, \\
\Kca_{\lambda_4}=\left\{ \left[\begin{smallmatrix} T(A)&\\& 1 \end{smallmatrix}\right], \left[\begin{smallmatrix}S(A)&\\& -1   \end{smallmatrix}\right]: \; A \in \U (2)  \right\},
\end{array}
$$
\noindent  where if for any $z=a+ib\in \CC$ we denote by
  $ M(z)=\left[\begin{smallmatrix} a&b\\ -b&a \end{smallmatrix}\right]$ and
  $\tilde{M}(z)=\left[\begin{smallmatrix} a&b\\ b&-a \end{smallmatrix}\right],$
  then for any $A=\left[\begin{smallmatrix}z_1 &z_2\\z_3&z_4 \end{smallmatrix}\right] \in \U(2),$ $T(A)$ and $S(A)$ are the matrices in $\glg(4)$ given by
$$
\begin{array}{lcl}
   T(A)=\left[\begin{smallmatrix} M(z_1)&M(z_2)\\ M(z_3)&M(z_4) \end{smallmatrix}\right], & \qquad
  S(A)=\left[\begin{smallmatrix}  \tilde{M}(z_1)&\tilde{M}(z_2)\\ \tilde{M}(z_3)&\tilde{M}(z_4)\end{smallmatrix}\right].
   \end{array}
$$
It is easy to see that by considering
$A_1=\left[\begin{smallmatrix}
  0&1\\1&0
\end{smallmatrix}\right]$ and
$A_2=\left[\begin{smallmatrix}
  i&0\\0&i
\end{smallmatrix}\right]$ then the element in $\Kca_{\lambda_4}$ corresponding to $T(A_1)$ permutes $a, b$ with $c, d$ and the one corresponding to $T(A_2)$ permutes $a$ with $b$ and fixed $c,d.$
Therefore $W_{\mu_4}$ is generated by $\{ \Id, \sigma_{12}, \sigma_{34}, \sigma_{13}\sigma_{24}\}.$

It can be seen that any derivation $\diag(a,b,c,a+b-c,a+b)$ is related under the action of $W_{\mu_4}$ to a multiple of an element in
$$
 \left\{\diag(a,b,c,a+b-c,a+b): a>0,\; 2 c \ge a+b \ge 0,\,\, a \ge c \ge b\right\}.
$$
Therefore, the set $\PP\ag_{\lambda_4}/W_{\lambda_4}$ of corresponding $6$-dimensional solsolitons is parameterized by:
\begin{equation}
  \mathfrak{F}^2_2=\left\{(a,b,c)\in S^2: a>0,\; 2 c \ge a+b \ge 0,\,\, a \ge c \ge b   \right\},
\end{equation}
\noindent where the correspondence is given by $(a,b,c) \to \diag(a,b,c,a+b-c,a+b)$.

\noindent Einstein condition: $a=b=c$.

\begin{remark} {\rm As $\rank(\lambda_4) = 3,$ the set $\mathfrak{F}^2_2$ also parameterizes $\Gr_2(\ag_{\lambda_4})/W_{\lambda_4},$ the space of $7$-dimensional solsolitons with nilradical $\lambda_4.$}
\end{remark}

\vs

$\bullet$  $\mathbf{\lambda}_5$ (free $2$-step nilpotent Lie algebra with $3$ generators).
$$\pg_{\lambda_5}=\left\{\diag(\left[\begin{smallmatrix} a&\alpha\\ \alpha&b \end{smallmatrix}\right],a+b,\left[\begin{smallmatrix} 2a+b &\alpha\\ \alpha&a+2b \end{smallmatrix}\right]),\,:\, a,b, \alpha \in \RR \right\},\qquad \kg_{\lambda_5}=\{0\}.$$
So we can choose
$$\ag_{\lambda_5}=\left\{ \diag(a,b,a+b,2a+b,a+2b),\,:\, a,b \in \RR \right\}.$$
It is easy to see that by acting with $A=\diag(\left[\begin{smallmatrix} 0&1\\1&0 \end{smallmatrix}\right], -1, \left[\begin{smallmatrix}0&-1\\-1&0 \end{smallmatrix}\right]) \in \Kca_{\lambda_5},$ we can permute $a$ with $b.$
\noindent The corresponding $6$-dimensional solsolitons are parameterized by $\mathfrak{F}^1_1$ (see (\ref{f1})).

\noindent Einstein condition: $a=b$.
\vspace{.3cm}

$\bullet$ $\mathbf{\lambda}_6.$
$$\pg_{\lambda_6}=\left\{ \diag(a,\left[\begin{smallmatrix}b&\beta\\\beta&c\end{smallmatrix}\right],
\left[\begin{smallmatrix}a+b&\beta\\ \beta&a+c\end{smallmatrix}\right]):\, a,b,c,\beta \in \RR \right\},$$
 $$\ag_{\lambda_6}=\left\{ \diag(a,b,c,a+b,a+c):a,b,c \in \RR \right\}, \Kca_{\lambda_6}=\left\{ \diag(\varepsilon,H,\varepsilon H), \, H \in \Or(2), \varepsilon = \pm 1 \right\}.$$
The corresponding $6$-dimensional solsolitons are parameterized by: $(c,b,a) \in \mathfrak{F}^2_1$ (see (\ref{f2})).

\noindent Einstein condition: $\frac{3}{2}\,a=b=c$.
\vspace{.3cm}

$\bullet$  $\mathbf{\lambda}_7$  ($\simeq \hg_3 \oplus \RR X_3 \oplus \RR X_4$).
$$\begin{array}{c}
  \pg_{\lambda_7}=\left\{ \diag(\left[\begin{smallmatrix}a&\alpha\\ \alpha& b\end{smallmatrix}\right],\left[\begin{smallmatrix}  c& \beta\\ \beta & d
  \end{smallmatrix}\right],a+b):\, a,b,c,d,\alpha, \beta \in \RR \right\},\\
\Kca_{\lambda_7}=\left\{ \diag(H_1, H_2, \det H_1), \,\, H_i \in \Or(2)  \right\},
\ag_{\lambda_7}=\left\{ \diag(a,b,c,d,a+b):\, a,b,c,d \in \RR \right\}.
\end{array}$$
$W_{\lambda_7}:$ group generated by $\{\Id, \sigma_{12},\sigma_{34}\}.$

\noindent The corresponding $6$-dimensional solitons are parameterized by:
\begin{equation}\label{f4}
  \mathfrak{F}_1^3=\left\{(a,b,c,d)\in S^3:   \begin{array}{l}
    a > 0,\, a \ge b, c \ge d, \, \mbox{or }\\
    a=b=0,\, c \ge |d|
  \end{array}\right\}.
\end{equation}
 Einstein condition: $a=b,$ $c=d,$ $3a=2c$.
\vspace{.3cm}

$\bullet$  $\mathbf{\lambda}_8$ ($\simeq\eta_1 \oplus \RR X_3$).
$$\begin{array}{l}
  \ag_{\lambda_8} = \pg_{\lambda_8}=\left\{ \diag(a,b,c,a+b,2a+b):\, a,b,c \in \RR \right\},\quad
\Kca_{\lambda_8}\subseteq \diag(\pm1, \dots, \pm1).
\end{array}$$
The corresponding $6$-dimensional solsolitons are parameterized by $\PP\RR^3$. Recall that $\PP\RR^3$ can be parameterized as
$$\PP\RR^3=\left\{(x_1,x_2,x_3) \in S^{2}: \begin{array}{l}
  x_3> 0,\, \mbox{or}\, \\  x_3 = 0,\, (x_1,x_2)\in \PP\RR^{2}
  \end{array}\right\},
$$
where $S^{n}$ is the $n$-dimensional sphere in $\RR^{n+1}$ and $\PP\RR^2$ is seen as in (\ref{pr2}) below.

\noindent Einstein condition: $2a=b=\frac{2}{3}c.$
\vs


\subsection{Nilsolitons of dimension $6$}(see Tables \ref{pc23} and \ref{pc45})

\vs

$\bullet$  $\mathbf{\mu}_6.$
$$\ag_{\mu_6}=\pg_{\mu_6}=\left\{\diag(a,b,a+b,2a+b,a+2b,3a+b),\; a,b \in \RR\right\},$$
$$\Kca_{\mu_6}\subseteq \diag(\pm1, \dots, \pm1).$$
The space of corresponding  $7$-dimensional solsolitons is parameterized by $\PP\RR^2,$ that can be seen as
\begin{equation}\label{pr2}
  \PP\RR^2=\{(x,y) \in S^1: y\ge 0,\, x > -1\}.
\end{equation}
\noindent Einstein condition: $2 a=b.$
\vs

\begin{remark} As we shall see, although  $\mu_7$ and $\mu_8$ are very similar, $\pg_{\mu_7}$ is completely different from $\pg_{\mu_8}$ (see Table \ref{sol6D1}). In particular, $\pg_{\mu_7}$ has non-diagonal elements (with respect to our fixed basis) but $\pg_{\mu_8}$ does not. The same holds for  $\mu_{15}$ and $\mu_{16}.$
\end{remark}

$\bullet$ $\mathbf{\mu}_7.$
$$\ag_{\mu_7} = \pg_{\mu_7}=\left\{\diag(\left[\begin{smallmatrix}a&b\\ b&a \end{smallmatrix}\right],2a,\left[\begin{smallmatrix}3a&b\\ b&3a \end{smallmatrix}\right],4a): a,b \in \RR\right\}. $$
$$\Kca_{\mu_7} = \left\{\diag(A,\varepsilon,\varepsilon A,\varepsilon): \varepsilon = \det A,
A= \left[\begin{smallmatrix}0&\pm 1\\ \pm 1&0 \end{smallmatrix}\right] \right\}.$$
The corresponding set of $7$-dimensional solsolitons is parameterized by $\PP\RR^2.$

\noindent Einstein condition: $b=0.$

\vs

$\bullet$ $\mathbf{\mu}_{15}.$ For simplicity we will change the basis to $\beta=\{X_1,X_3,X_2,X_4,X_5,X_6\}.$ We then get that
$$\ag_{\mu_{15}} = \pg_{\mu_{15}}=\left\{\diag(\left[\begin{smallmatrix}a&c\\c&a \end{smallmatrix}\right],b,\left[\begin{smallmatrix}a+b&-c\\-c&a+b \end{smallmatrix}\right], 2a+b): a,b,c \in \RR\right\}. $$
It is not hard to see that the action of $\Kca_{\mu_{15}}$ can at most permute $c$ with $-c$ (that is $a+c$ with $a-c$) and this can be done by acting with
$A=\diag(\left[\begin{smallmatrix}0&1\\-1&0 \end{smallmatrix}\right],1,\left[\begin{smallmatrix}0&-1\\1&0 \end{smallmatrix}\right],-1) \in \Kca_{\mu_{15}}.$

\noindent The corresponding set of  $7$-dimensional solsolitons is parameterized by:
\begin{equation}
  \mathfrak{F}^2_3=\left\{(a,b,c) \in S^2: \begin{array}{ll}
    a>0,\, c\ge 0,\, \mbox{or} \\ a=0,\; b,c\ge 0
  \end{array}\right\}.
    \end{equation}
Einstein condition: $a=b, c=0.$
\vs

$\bullet$ $\mathbf{\mu}_{16}.$ As in the previous case, we will change the basis to $\beta$ to get nicer matrices. In this way we get that
$$\ag_{\mu_{16}} = \pg_{\mu_{16}}=\left\{\diag(a,a,b,a+b,a+b,2a+b): a,b \in \RR\right\}, $$
$$\Kca_{\mu_{16}}=\left\{\diag(H_1,\pm 1,H_2, \pm 1): H_i \in \Or(2)\right\}.$$
The corresponding set of $7$-dimensional solsolitons is parameterized by $\PP\RR^2.$

\noindent Einstein condition: $a=b.$

\vs

$\bullet$ $\mathbf{\mu}_{17}$  ($\simeq\lambda_5 \oplus \RR X_3).$
$$   \pg_{\mu_{17}}=\left\{\diag(\left[\begin{smallmatrix}a&d\\d&b \end{smallmatrix}\right],c,a+b,\left[\begin{smallmatrix}2a+b&d\\d&a+2b \end{smallmatrix}\right]): a,b,c,d \in \RR\right\},$$
$$\ag_{\mu_{17}} = \left\{\diag(a,b,c,a+b,2a+b,a+2b): a,b,c \in \RR\right\},$$
$$\Kca_{\mu_{17}}=\left\{ \diag(A,\pm 1,\varepsilon,\varepsilon A): A \in \Or(2),\varepsilon=\det A \right\}.$$
The corresponding set of $7$-dimensional solsolitons is parameterized by $\mathfrak{F}^2_1$ (see (\ref{f2})).

 \noindent Einstein condition $a=b=\frac{5}{12}c.$


$\bullet$ $\mathbf{\mu}_{20}.$

$\qquad \qquad \ag_{\mu_{20}} = \pg_{\mu_{20}}=\left\{\diag(a,b,a+b,a+b,2a+b,a+2b):\; a,b \in \RR\right\}.$

 A straightforward calculation shows that one can permute $a$ with $b$ by acting with \\ $A=\diag(\left[\begin{smallmatrix}&1\\1& \end{smallmatrix}\right],\left[\begin{smallmatrix}\frac{1}{2}&\frac{\sqrt{3}}{2}\\\frac{\sqrt{3}}{2}&-\frac{1}{2} \end{smallmatrix}\right],\left[\begin{smallmatrix}&-1\\-1& \end{smallmatrix}\right]) \in \Kca_{\mu_{20}}.$

    \noindent The corresponding set of $7$-dimensional solsolitons is parameterized by $\mathfrak{F}^1_1$ (see (\ref{f1})).

    \noindent Einstein condition: $a=b.$

\vs

$\bullet$  $\mathbf{\mu}_{24}$  (free $3$-step nilpotent Lie algebra with $2$ generators).
$$\pg_{\mu_{24}}=\left\{\diag(\left[\begin{smallmatrix}a&b&c\\b&d&e\\c&e&f \end{smallmatrix}\right],\left[\begin{smallmatrix} a+d&e&-c\\e&a+f&b\\-c&b&d+f \end{smallmatrix}\right]):\;a,b,c,d,e,f \in \RR\right\},$$
$$\ag_{\mu_{24}}=\left\{ D(a,b,c)= \diag(a,b,c,a+b,a+c,b+c): a,b,c \in \RR\right\}.$$
It is known that $\Kca_{\mu_{24}} \simeq\Or(3)$ and from this we can see that the corresponding set of $7$-dimensional solsolitons is parameterized by $\PP \RR^3/ S_3$. This set can be explicitly realized as
$$\mathfrak{F}_4^2=\left\{(a,b,c)\in S^2: \begin{array}{l}
    b>0, a\ge b \ge c, \mbox{or} \\
    b=0,\; |c|\le a
  \end{array}\right\}.$$
\noindent (see also \cite{GP}). It will be Einstein if and only if $a=b=c.$

\vs

$\bullet$  $\mathbf{\mu}_{25}.$
 $$ \pg_{\mu_{25}}= \{\diag(a,b,\left[\begin{smallmatrix}c&\beta\\\beta&2a+b-c \end{smallmatrix}\right],a+b,2a+b): a,b,c, \beta \in \RR\},$$
$$\ag_{\mu_{25}}=\left\{ \diag(a,b,c,2a+b-c,a+b,2a+b):\, a,b,c \in \RR\right\}, $$
$$\Kca_{\mu_{25}}= \left\{ \diag(\nu,\varepsilon,H, \varepsilon\nu,\varepsilon), H\in \Or(2), \varepsilon =\det H, \nu=\pm 1 \right\}.$$
The corresponding set of $7$-dimensional solsolitons is parameterized by $\PP\RR^3$.

\noindent Einstein condition: $a/5=b/8=c/9$

\vs

$\bullet$  $\mathbf{\mu}_{26}$ ($\simeq \eta_1\oplus (\RR X_3 \oplus \RR X_4)).$
 $$ \pg_{\mu_{26}}= \{\diag(a,b,C,d,e): a,b,d,e \in \RR, C \in \sym(2)\},$$
$$\ag_{\mu_{26}}=\left\{\diag(a,b,c,d,a+b,2a+b):\, a,b,c,d \in \RR\right\},$$
$$\Kca_{\mu_{26}}= \left\{ \diag(\varepsilon,\nu,H,\varepsilon\nu,\nu), H\in \Or(2),\epsilon=\pm 1,\nu=\pm 1 \right\}.$$
$W_{\mu_{26}}$: group generated by $\{\Id, \sigma_{34}\}.$

\noindent The corresponding set of $7$-dimensional solsolitons is parameterized by (see (\ref{f2})):
\begin{equation}\label{f5}
  \mathfrak{F}_2^3=\left\{(a,b,c,d)\in S^3: \begin{array}{l}
    a>0,\, c\le d, \mbox{or} \\
    a=0,\,(c,d,b) \in \mathfrak{F}_1^2
  \end{array}\right\}.
\end{equation}
Einstein condition: $6a=3b=2c=2d.$

\vs

$\bullet$ $\mathbf{\mu}_{27}$ ($\simeq \lambda_3\oplus \RR X_3$).
$$\ag= \pg_\mu(\mu_{27})= \left\{\diag(a,b,c,2b,a+b,a+2b):\, a,b,c \in \RR\right\},\,
\Kca(\mu_{27}) \subset \diag(\pm 1,\dots, \pm 1).$$
The corresponding set of $7$-dimensional solsolitons is parameterized by $\PP\RR^3.$

\noindent Einstein condition: $4a= 3b=2c.$
\vs

$\bullet$ $\mathbf{\mu}_{28}.$
$$ \pg_{\mu_{28}}= \left\{ \diag(\left[\begin{smallmatrix}\begin{smallmatrix}a&\\&a\end{smallmatrix}&B\\B^t& \begin{smallmatrix}b&\\&b\end{smallmatrix} \end{smallmatrix}\right],a+b,a+b):
B=\left[\begin{smallmatrix}\alpha &\beta\\-\beta&\alpha \end{smallmatrix}\right], a,b, \alpha,\beta \in \RR,\right\}$$
$$\ag_{\mu_{28}}=\{\diag(a,a,b,b,a+b,a+b): a,b \in \RR\}$$
The group $\Aut(\mu_{28})$ has been calculated in \cite{S}, although it is easy to see that we can interchange $a$ with $b$ (the only possible permutation) by conjugating by $A=\left[\begin{smallmatrix}\begin{smallmatrix} &&&1\\&&1&\\&1&&&\\1&&& \end{smallmatrix}& \\ & \begin{smallmatrix}1&\\&-1\end{smallmatrix}\end{smallmatrix}\right] \in \Kca_{\mu_{28}}.$

\noindent The corresponding set of $7$-dimensional solsolitons is parameterized by $\mathfrak{F}^1_1$ (see (\ref{f1})).

\noindent Einstein condition: $a=b.$

\vs

$\bullet$ $\mathbf{\mu}_{29}.$
$$ \pg_{\mu_{29}}= \left\{ \diag(\left[\begin{smallmatrix} a &\alpha\\ \alpha &b \end{smallmatrix}\right], \left[\begin{smallmatrix}c&-\alpha\\ -\alpha&b+c-a \end{smallmatrix}\right],a+b,b+c):
a,b,c,\alpha \in \RR\right\},$$
$$\ag_{\mu_{29}}= \left\{\diag(a,b,c,c-a+b,a+b,b+c):
a,b,c\in \RR\right\},$$
$$\Kca_{\mu_{29}}= \left\{\diag(\left[\begin{smallmatrix}a&b\\c&d \end{smallmatrix}\right], \left[\begin{smallmatrix}d&c\\b&a \end{smallmatrix}\right],\varepsilon,1): H=\left[\begin{smallmatrix}a&b\\c&d \end{smallmatrix}\right]\in \Or(2), \varepsilon =\det H \right\}.$$
The corresponding set of $7$-dimensional solsolitons is parameterized by $\PP\RR^3$.

\noindent Einstein condition: $a=b,\, 4a=3c.$

\vs

$\bullet$ $\mathbf{\mu}_{30}$ ($\simeq \hg_3\oplus \hg_3$).
$$\pg_{\mu_{30}}= \left\{ \diag(\left[\begin{smallmatrix} a &\alpha\\ \alpha & b \end{smallmatrix}\right], \left[\begin{smallmatrix}c&\beta\\
\beta& d \end{smallmatrix}\right],a+b,c+d):
a,b,c,d,\alpha,\beta \in \RR\right\},$$
$$\ag_{\mu_{30}}= \left\{\diag(a,b,c,d,a+b,c+d):
a,b,c,d,\in \RR\right\},$$
$$\Kca_{\mu_{30}}^0= \left\{\diag(H_1,H_2,\varepsilon_1,\varepsilon_2): \; H_i \in\Or(2), \varepsilon_i =\det H_i \;i=1,2 \right\}.$$
In this case there is another orthogonal automorphism that act non trivially on subspaces of $\ag$ given by
$H=\diag( \left[\begin{smallmatrix} 0&H'\\H'&0\end{smallmatrix}\right],H')$ where $H'= \left[\begin{smallmatrix} 0&1\\1&0\end{smallmatrix}\right]$, and therefore we get that $W_{\mu_{30}}$ is the group generated by $\{\Id, \sigma_{12},\sigma_{34},\sigma_{13}\sigma_{24}\}.$ It is easy to see that $\RR^4/ W_{\mu_{30}}$ can be parameterized by
$$\mathfrak{D}=\left\{(a,b,c,d)\in \RR^4: a\ge b, c\ge d, \begin{array}{l} a+b > c+d \mbox{ or }\\ a+b=c+d \mbox{ and } a \ge c\end{array}
\right\}.$$
 Hence, the corresponding set of $7$-dimensional solsolitons is parameterized by
 $$\mathfrak{F}_3^3=\left\{(a,b,c,d)\in S^3 \cap \mathfrak{D}: a>0,  a+b \ge 0, \mbox{ if } d,\,c+d < 0 \mbox{ then } a > |d|
\right\}.$$
Einstein condition: $a=b=c=d$.

\vs

$\bullet$ $\mathbf{\mu}_{31}$ ($\simeq \lambda_6\oplus \RR X_4).$
$$ \pg_{\mu_{31}}= \left\{ \diag(a,\left[\begin{smallmatrix}b&\alpha\\\alpha&c \end{smallmatrix}\right], d, \left[\begin{smallmatrix}a+b&\alpha\\ \alpha&a+c \end{smallmatrix}\right]):
a,b,c,d,\alpha \in \RR\right\},$$
$$\ag_{\mu_{31}}= \left\{\diag(a,b,c,d,a+b,a+c): \;a,b,c,d,\in \RR\right\},$$
$$\Kca_{\mu_{31}}= \left\{ \diag(\varepsilon,H,\pm 1,\varepsilon H), H\in \Or(2), \varepsilon=\pm 1\right\}.$$

\noindent $W_{\mu_{31}}:$  generated by $\{\Id, \sigma _{23}\}.$

\noindent The corresponding set of $7$-dimensional solsolitons is parameterized by: $(a,d,b,c) \in \mathfrak{F}_2^3$ (see (\ref{f5})).

\noindent Einstein condition: $4b=4c=6a=3d.$

\vs

$\bullet$ $\mathbf{\mu}_{32}$ $(\simeq \lambda_4\oplus \RR X_5$) (see $\lambda_4$).
$$ \pg_{\mu_{32}}= \left\{\diag(\left[
\begin{smallmatrix}
\begin{smallmatrix}a&\alpha\\ \alpha & b \end{smallmatrix} & B\\ B &
\begin{smallmatrix} c &\beta\\ \beta & a+b-c\end{smallmatrix}\end{smallmatrix}\right],d,a+b):\, a,b,c,d,\alpha,\beta \in \RR, \; B= \left[\begin{smallmatrix}u&v\\v & -u \end{smallmatrix} \right] \right\},$$
$$ \ag_{\mu_{32}}= \left\{\diag(a,b,c,a+b-c,d,a+b):\,
a,b,c,d,\in \RR\right\}.$$
Concerning $\Kca_{\mu_{32}}$ we can use the information we have from $\lambda_4$ to obtain that
the corresponding set of $7$-dimensional solsolitons is parameterized by:
$$\mathfrak{F}_4^3=\left\{(a,b,c,d) \in S^3: \begin{array}{l} d>0,\, a \ge c \ge b,\; 2c \ge a+b, \;\; \mbox{or} \\  d=0,\, (a,b,c) \in \mathfrak{F}_2^2 \end{array} \right\}.$$
Einstein condition: $a=b=c=\frac{3}{4}d.$

\vs

$\bullet$ $\mathbf{\mu}_{33}$ ($\simeq \hg_3\oplus (\RR X_3 + \RR X_4 + \RR X_5)).$
$$ \pg_{\mu_{33}}= \left\{ \diag(A,B,\tr A): A\in \sym(2), B \in \sym(3) \right\},$$
$$\ag_{\mu_{33}}= \left\{\diag(a,b,c,d,e,a+b):\, a,b,c,d,e\in \RR\right\},$$
$$\Kca_{\mu_{33}}= \left\{\diag(H_1,H_2,\varepsilon):\, H_1\in \Or(2),H_2\in \Or(3), \varepsilon =\det H_1\right\}.$$
The corresponding set of $7$-dimensional solsolitons is parameterized by (see \ref{f4}):
$$\mathfrak{F}^4_{1}=\left\{(a,b,c,d,e) \in S^4: \begin{array}{l} e>0,\, a\ge b,\, c \le d \le e, \, \mbox{or} \\ e=0,\, (a,b,c,d) \in \mathfrak{F}_1^3  \end{array} \right\}.$$
Einstein condition: $3a=3b=2c=2d=2e.$


\vs

\section{The space of $m$-dimensional solsolitons}

   We now describe the space of solsolitons for each dimension $\le 6,$ by using the information we have obtained in the previous sections.
 Let $\ngo_\mu=(\RR^n,\mu)$ denote a metric nilpotent Lie algebra as in Sections $3$ and $4$ (recall that we have fixed a basis and an inner product on $\RR^n$). Hence we may denote as
$${\mathcal Nil}(n)= \{\ngo_\mu=(\RR^n,\mu):\; \mu \mbox{ is a nilsoliton} \}$$
the moduli space of nilsolitons of dimension $n$ up to isometry and scaling.

According to Proposition \ref{const} and (\ref{para}) we have that ${\mathcal Sol}(m)$, the moduli space of $m$-dimensional solsolitons (up to isometry and scaling), is given by:
$$
{\mathcal Sol}(m)= \bigcup_{\begin{array}{c}\mu \in {\mathcal Nil}(n)\\ m\ge n \ge m/2\end{array}}\;\; \big\{\ag \oplus \ngo_\mu:\;\; \ag \; \in\; \Gr_{m-n}(\ag_\mu)/{W_{\mu}}\big\},
$$
where $\Gr_r(\ag_\mu)$ is the Grassmanian of $r$-dimensional subspaces of $\ag_\mu$ (a maximal abelian subspace of symmetric derivations of $\ngo_\mu$), and $\ag \oplus \ngo_\mu$ is the (metric) solvable Lie algebra constructed in Proposition \ref{const}. Recall that $m-n \le k_\mu = \rank(\mu) = \dim(\ag_\mu).$

We note that each $\mu \in {\mathcal Nil}(m)$ is a point in ${\mathcal Sol}(m)$ and the same holds when $m=k_\mu + \dim \ngo_\mu,$ that is, there is a single point $\ag_\mu \oplus \ngo_\mu$ in the space of solsolitons of dimension $k_\mu+ \dim \ngo_\mu$ which is always Einstein (see Section \ref{maxrank}).

More explicitly, we can see this moduli space as
$$ \begin{array}{l}
{\mathcal Sol}(m)= {\mathcal Nil}(m) \;\;\cup  \displaystyle{\bigcup_{\mu \in {\mathcal Nil}(m-1)}}  \Gr_1(\ag_\mu)/{W_{\mu}}
\;\;\; \cup \displaystyle{\bigcup_{\tiny{\begin{array}{c}\mu \in {\mathcal Nil}(m-2) \\ k_\mu \ge 2 \end{array}}}}  \Gr_2(\ag_\mu)/{W_{\mu}}\\
\cup \dots
\cup  \displaystyle{\bigcup_{\tiny{\begin{array}{c}\mu \in {\mathcal Nil}(m-r)\\ k_\mu \ge r\end{array}}}}  \Gr_r(\ag_\mu)/{W_{\mu}}
\;\;\;\cup \dots \cup  \displaystyle{\bigcup_{\tiny{\begin{array}{c}\mu \in {\mathcal Nil}(m/2)\\ k_\mu = m/2\end{array}}}}  \Gr_r(\ag_\mu)/{W_{\mu}}
\end{array}
$$
for $m$ even. For $m$ odd, the same expression holds but without the last union.

We therefore get that:
$$\begin{array}{l}
{\mathcal Sol}(2)=  \left\{ \RR^2\right\} \cup \left\{(\ag_\RR \oplus \RR)\right\},\\ \\
{\mathcal Sol}(3) = \left\{\hg_3, \RR^3 \right\} \; \cup \;  \mathfrak{F}_1^1(\RR^2),\\
\\
{\mathcal Sol}(4) = \left\{\eta_1, \eta_2, \RR^4\right\}  \;\cup \;  \mathfrak{F}_1^1(\hg_3) \;\cup \; (\PP\RR^3/S_3)(\RR^3) \;\cup \;
\left\{\ag_{\RR^2} \oplus \RR^2\right\},\\
\\
{\mathcal Sol}(5) = \left\{\lambda_1, \dots,\lambda_8,\RR^5 \right\} \cup \PP\RR^2(\eta_1) \cup  \mathfrak{F}_1^2(\eta_2) \cup (\PP\RR^4/S_4)(\RR^4) \cup  \{\ag_{\hg_3} \oplus \hg_3\} \cup \mathfrak{F}_4^2(\RR^3),\\
\\
{\mathcal Sol}(6) = \left\{\mu_1, \dots, \mu_{33}, \RR^6 \right\}\, \cup \, \PP\RR^2(\lambda_1) \,\cup \, \PP\RR(\lambda_2) \,\cup\, \PP\RR^2(\lambda_3)\, \cup\, \mathfrak{F}_2^2(\lambda_4)\, \cup \,\mathfrak{F}_1^1(\lambda_5) \\ \\ \qquad \cup\, \mathfrak{F}_1^2(\lambda_6)\, \cup \,\mathfrak{F}_1^3(\lambda_7) \,\cup \, \PP\RR^3(\lambda_8) \cup  (\PP\RR^5/S_5)(\RR^5)
\cup \, \{\ag_{\eta_1} \oplus \eta_1\} \cup \, \mathfrak{F}_1^2(\eta_2)\\ \\ \qquad \cup \, (\Gr_2(\RR^4)/S_4)(\RR^4) \,\cup \, \{\ag_{\RR^3} \oplus \RR^3\},
\end{array}$$
where to distinguish among them we have added the metric to the notation. In this way, for example, $\PP\RR^2(\lambda_1)$ denotes the set $\PP\RR^2$ that parameterize the solsolitons associated to $\lambda_1$ in the corresponding dimension. We note that $\mathfrak{F}_i^j(\mu)$ is a real semialgebraic set and $\dim \mathfrak{F}_i^j(\mu) = j.$

Moreover, with the results we have obtained in Sections $3$ and $4$, one can describe ${\mathcal Sol}(7).$ It is given by the union of the following:
\begin{itemize}
  \item one point for each element of ${\mathcal Nil}(7);$
  \item six points coming from the $6$-dimensional nilsolitons from Table \ref{sol6D1}, eight copies of $\PP\RR^2$ and three copies of $\PP\RR^3$ coming from the $6$-dimensional algebras in Table \ref{sol111}, and the one-dimensional extensions of the elements in ${\mathcal Nil}(6),$ described in Section $4$;
\item three points corresponding to $\lambda_1$,$\lambda_3$ and $\lambda_5$, $ \mathfrak{F}_2^2(\lambda_4),\, \mathfrak{F}_1^2(\lambda_6),\, \PP\RR^3(\lambda_8),\\ (\Gr_2(\RR^4)/{W_{\lambda_7}})(\lambda_7)$ and $(\Gr_2(\RR^5)/S_5)(\RR^5);$
\item one point corresponding to $\eta_2$ and $(\Gr_3(\RR^4)/S_4)(\RR^4).$
\end{itemize}

 We note that in this union not everything is explicitly described. Indeed, on the one hand we have ${\mathcal Nil}(7)$  which, as far as we know, has not been completely classified yet. On the other hand, for $r>0, $ besides from the abelian cases, there is just one set which is not explicitly described in this paper: $(\Gr_2(\RR^4)/W_{\lambda_7})(\lambda_7).$

\section{On negatively curved solsolitons.}

To conclude we will make some remarks on the curvature of solsolitons.
Recall that for each nilsoliton $\mu \in {\mathcal Nil}(n)$ and each $1 \le r < k_\mu = \rank(\mu),$ we have a family $X_{\mu,r}$ of solsolitons obtained from $\mu$ as in Proposition \ref{const}. Moreover, this family is parameterized, up to isometry and scaling, by
$$ X_{\mu,r} = \Gr_r(\ag_\mu)/W_\mu.$$
We then have that $X_{\mu,r} \subset {\mathcal Sol}(m)$ where $m=r+n,$ and since $W_\mu$ is a finite group, this quotient inherits many of the properties of $\Gr_r(\ag_\mu).$ In particular, $X_{\mu,r}$ is always connected and compact. Recall that if $\rank(\mu)=1,$ there is only one solsoliton associated to $\mu$ (up to isometry and scaling) and therefore $X_{\mu,1}$ is just a single point, as it always happens for $X_{\mu,k_\mu}$.

\begin{remark}
Since a solsoliton $\sg=\ag\oplus \ngo_\mu \in X_{\mu,r}$ is Einstein if and only if $D_1 \in \ag$ (see Proposition \ref{const}),  then for $1 \le r < k_\mu$ the set of non-Einstein solsolitons is open and dense in $X_{\mu,r}.$
\end{remark}

\begin{remark}
Let $\sg_0 \in X_{\mu,r}$ be an Einstein solvmanifold with sectional curvature $K_{\sg_0}<0.$  By continuity, if $\sg \in X_{\mu,r}$ is sufficiently closed to $\sg_0$ then $K_\sg <0$ as well. The classification of Einstein solvmanifolds with non-positive sectional curvature of dimension $\le 6$ is given in \cite{NN}. It is proved there that the Einstein solvmanifolds with negative sectional curvature of dimension $\le 6$ are the symmetric spaces $\RR H^6$ and $\CC H^6$ and the Einstein one-dimensional extensions of $\lambda_6$ and $\lambda_7$ (see \cite[Theorem 2]{NN} and Table \ref{5dn}). Let us look at $\lambda_6$, for example. We have that $\rank(\lambda_6)=3$ and $X_{\lambda_6,1} \simeq \mathfrak{F}_1^2$ (see \ref{f2}), and therefore around $\sg_0=\RR D_1 \oplus \ngo_{\lambda_6} \in X_{\lambda_6,r}$ we get a $2$-dimensional family of $6$-dimensional non-Einstein solsolitons with $K_{\sg}<0$.
\end{remark}

\begin{remark}
 We can apply the same argument to the Ricci curvature. In fact, if $\sg_0 \in X_{\mu,r}$ is Einstein then $\ricci_{\sg_0}=c \Id <0$ and therefore  $\ricci_{\sg}<0$ for any $\sg \in X_{\mu,r}$ sufficiently close to $\sg_0.$

On the other hand, if $A \in \ag_\mu$ is orthogonal to $D_1$ then $\sg_A = \RR A \oplus \ngo_\mu \in X_{\mu,1}$ does not have negative Ricci curvature. Indeed, it is easy to see that $\langle D_1, A \rangle =0$ implies that ${\ricci_{\sg_A}}|_{\ngo_\mu}= \ricci_\mu$ and therefore $\ricci_{\sg_A}$ has both signs. Note that therefore  $\sg_A$ has not negative sectional curvature either. The same holds if $\ag \subset \ag_\mu$ is a subspace orthogonal to $D_1.$ For the existence of these orthogonal subspaces we need $\rank(\mu) >1$.

Summarizing, if $1 \le r < \rank(\mu)$ then there are always some solsolitons in  $X_{\mu,r}$ with negative Ricci curvature and others with positive and negative directions for the Ricci tensor.
\end{remark}

\begin{remark}
For $\mu=0$ and $r=1,$ the Einstein solsoliton corresponds to
the real hyperbolic space $\RR H^{n+1},$ which is a symmetric space of constant negative sectional curvature. According to the above observations, in $X_{0,1}= \PP \RR^n/S_n,$ there are families depending on $n-1$ parameters of non-Einstein solsolitons with $K_\sg <0$ as closed as one wishes to $\RR H^{n+1}.$

\end{remark}


\end{document}